\documentclass[runningheads,a4paper]{llncs}
\usepackage{epsfig,color}
\usepackage{amsmath,amssymb}

\usepackage{graphicx}

\usepackage{url}
\urldef{\mailsa}\path|{alfred.hofmann, ursula.barth, ingrid.haas, frank.holzwarth,|
\urldef{\mailsb}\path|anna.kramer, leonie.kunz, christine.reiss, nicole.sator,|
\urldef{\mailsc}\path|erika.siebert-cole, peter.strasser, lncs}@springer.com|

   \def\1{{\mathbbm 1}}

   \def\R{{\mathbb R}}

   \def\Pf{{\it Proof.$\;$}}

   \def\qed{\hspace{10cm}$\diamond$}

   \def\lin{\mbox{\rm lin}}

   \def\({\langle}
   \def\){\rangle}
   \def\mb{\boldsymbol}

   \def\cM{{\mathcal M}}

   \def\cP{{\mathcal P}}

   \def\cV{{\mathcal V}}

   \def\1{\mb1}

   \def\v0{{\bf 0}}

   \def\ov{\overline}

\begin{document}

\mainmatter

\title{Max-Flow on Regular Spaces}

\author{
Ulrich Faigle\inst{1}
\and Walter Kern\inst{2}
\and
Britta Peis\inst{3}
}

\authorrunning{Regular max-flow}

\institute{Math. Institut, Universit\"at zu K\"oln,
Weyertal 80, D-50931 K\"oln, \email{faigle@zpr.uni-koeln.de} \and
Universiteit Twente, P.O. Box 217, NL-7500 AE Enschede, \\\email{w.kern@math.utwente.nl} \and
Technische Universit\"at Berlin, Stra\ss{e} des 17.\ Juni 135, D-10623 Berlin,
\email{peis@math.tu-berlin.de}
}

\maketitle
\abstract{The max-flow and max-coflow problem on directed graphs is studied in the common generalization to regular spaces, \emph{i.e}, to kernels or row spaces of totally unimodular matrices. Exhibiting a submodular structure of the family of paths within this model we generalize the Edmonds-Karp variant of the classical Ford-Fulkerson method and show that the number of augmentations is quadratically bounded if augmentations are chosen along shortest possible augmenting paths.}

\section{Introduction}\label{sec:Introduction}
Let $G=(V,E)$ be a (finite) oriented graph with vertex set $V$, arc set $E$ and a distinguished return arc $r=(t,s)\in E$ from the ''sink'' $t\in V$ to the ''source'' $s\in V$. Let $A\in \{-1,0,+1\}^{V\times E}$ be the vertex-arc incidence matrix of $G$. Referring to the vectors $f\in \ker A$ as \emph{flows} (a.k.a. \emph{circulations}) on $G$, the classical max-flow problem of Ford and Fulkerson~\cite{FF56} asks for a feasible flow $f$ of maximal value $f_r$ on $r$,
$$
    \max f_r \quad\mbox{s.t.}\quad f\in \ker A, 0\leq f\leq c,
$$
where $c:E\setminus\{r\}\to\R_+$ (and $c_r=\infty$) describes a capacity restriction on the arc set $E$. It is well-known that the augmenting path algorithm of Ford and Fulkerson can be made polynomial {\it via} the Edmonds-Karp variant that always augments along a shortest possible $(s,t)$-path (see, \emph{e.g.}, \cite{AMO93}).

\medskip
Dually, we may consider the row space $\lin~A$ of $A$, whose elements are the \emph{co-flows} (a.k.a. \emph{tensions}) on $G$ that are induced by vertex potentials $p:V\to\R$, and investigate the corresponding \emph{co-flow problem}
$$
    \max f_r \quad\mbox{s.t.}\quad f\in \lin A,\; 0\leq f\leq c.
$$

The latter problem seems to have received much less attention in the literature, in  spite of the fact that already Minty~\cite{Minty66} proposed a common generalization of the flow and co-flow problem on directed graphs to regular matroids by observing the total unimodularity of the vertex-arc incidence matrix $A$ and formulating the problem with respect to general totally unimodular matrices (see also Hoffman~\cite{Hoffman-abstract-flows} for a further abstraction).
Referring to a vector space $\cV\subseteq \R^n$ as \emph{regular} if it is the kernel (or row space) of some totally unimodular matrix, we study the (co-)flow problem in regular spaces (\emph{i.e.}, in Minty's generalized model).

\medskip
We show that the collection $\cP$ of ''paths'' in that general context has a natural algebraic structure that is ''submodular'' with respect to path length and allows us to ''uncross'' nonconformal paths. Taking advantage of this structure, we prove that the number of augmentations in the general model is bounded by $|E|^2$ if augmentation is always carried out along a shortest possible path, thus establishing a polynomial version of a Ford-Fulkerson type  approach to co-flows.

\section{Regular spaces}\label{sec:Regular-spaces}
We review briefly well-known properties of regular spaces (see, \emph{e.g.}, Tutte~\cite{Tutte71} (or \cite{BK92}) for more details). Let us fix a real vector space $\cV\subseteq \R^n$ and the corresponding index ground set $E =\{1,\ldots,n\}$. Recall that the \emph{support} of a vector $x\in \R^n$ is the set
$$
\|x\|:= \{j\in E\mid x_j\neq 0\}.
$$
A nonzero $x\in \cV$ is called \emph{elementary} if there is no nonzero $y\in \cV$ whose support is a proper subset of $\|x\|$. A \emph{primitive} vector is an elementary vector $p\in \cV$ with components $p_j\in \{-1,0,+1\}$. Elementary vectors are determined by their supports (up to scaling factors):

\medskip
\begin{lemma}\label{l.elementary} Let $x,y\in \cV$ be elementary with $\|x\|=\|y\|$. Then $y =\lambda x$ holds for some $\lambda \in \R$.

\qed
\end{lemma}

\medskip
$\cV$ is said to be \emph{regular} if each elementary vector is a scalar multiple of a primitive vector. Moreover, the following statements are equivalent:
\begin{enumerate}
\item[(i)] $\cV$ is a regular space.
\item[(ii)] $\cV$ is the row space of some totally unimodular matrix.
\item[(iii)] $\cV$ is the kernel of some totally unimodular matrix.
\end{enumerate}

\medskip
Let $x,y\in \R^n$ be arbitrary with components $x_j$ and $y_j$. $x$ is said to \emph{conform} to $y$ if
$$
     x_j y_j >0 \quad\mbox{whenever $x_j\neq 0$.}
$$

\medskip
For our purposes the following properties of regular spaces are relevant.

\medskip
\begin{lemma}\label{lemma1} Let $\cV$ be a regular space. Then  any $x\in \R^n$ is the sum of elementary vectors, each conforming to $x$.

\qed
\end{lemma}

\medskip
\begin{lemma}\label{lemma2} Let $\cV$ be a regular space and $x\in \R^n$. Then $x$ is integral (\emph{i.e.}, all components $x_j$ of $x$ are integers) if and only if $x$ is a sum of primitive vectors, each conforming to $x$.

\qed
\end{lemma}

\section{The regular max-flow problem}\label{sec:max-flow}
Let $\cV\subseteq \R^n$ be a regular space and refer to any element $f\in \cV$ as a \emph{flow} on $E=\{1,\ldots,n\}$. Let us fix an element $r\in E$ and a \emph{capacity} $c:E\setminus\{r\}\to \R_+$ and say that  the flow $f$ is \emph{feasible} if
$$
      0\leq f_j \leq c_j \quad\mbox{for all $j\neq r$.}
$$

The corresponding \emph{max-flow problem} is the optimization problem
\begin{equation}\label{eq.flow-problem}
\max~ \{f_r\mid \mbox{$f$ is a feasible flow}\}.
\end{equation}

\medskip
We will approach problem (\ref{eq.flow-problem}) in the spirit of Ford-Fulkerson. So  we call a primitive vector $P\in \cV$ an \emph{$r$-path} (or simply \emph{path} if the reference to $r$ is clear) if $P_r = +1$. Let $\cP$ denote the collection of all $r$-paths. Given a feasible flow $f$, we say that $P\in \cP$ is \emph{augmenting}
if there is some $\varepsilon >0$ such that the flow $f'=f+\varepsilon P$ is feasible (and hence improves $f$ in view of $f'_r = f_r +\varepsilon$).

\medskip
\begin{theorem}\label{t.optimality} Let $\cV$ be a regular space on $E$ and $c:E\setminus\{r\}\to \R_+$ a capacity. Then the feasible flow $f$ is optimal for (\ref{eq.flow-problem}) if and only if $f$ does not admit an augmenting path $P\in \cP$.
\end{theorem}

{\small\medskip\Pf Clearly, an optimal $f$ cannot admit an augmenting path.
So assume that $f$ is not optimal and consider a feasible flow $f'\in \cV$ with $f_r'> f_r$.
By Lemma~\ref{lemma1}, the difference vector $y=f'-f$ is a sum of elementary vectors $x$,
each conforming to $y$. So $f+x$ is feasible for each of these elementary vectors.

\medskip
Because of $y_r >0$, at least one of them, say $\ov{x}$ must have $\ov{x}_r >0$.
Since elementary vectors are scalar multiples of primitive vectors, we have $\ov{x}=\lambda P$ for some $\lambda >0$ and $P\in \cP$. So $f+\ov{x}= f+\lambda P$ is feasible and exhibits $P$ as augmenting.

\qed}

\medskip
In view of Theorem~\ref{t.optimality}, the \emph{Ford-Fulkerson method} for the max-flow problem in a directed graph $G=(V,E)$ with a distinguished arc $r=(t,s)$ carries over to the regular max-flow problem:

\medskip
\begin{enumerate}
\item Start with the zero flow $f=0$.
\item If $f$ is the current feasible flow, determine an augmenting path $P\in \cP$ and $\varepsilon >0$ maximal so that $f'= f+ \varepsilon P$ is a feasible flow. Replace $f$ by $f'$ and repeat until a flow is attained that cannot be augmented.
\item Return the final flow $f$.
\end{enumerate}

\medskip
In the remainder of this note, we will show that the Ford-Fulkerson method can be carried with at most $|E|^2$ augmentations. In particular, we will show that the Edmonds-Karp technique of \emph{shortest} path augmentations of flows in directed graphs generalizes to regular spaces.

\section{Shortest path augmentation}\label{sec:shortest-path}
We evaluate the \emph{length} $|x|$ of a vector $x\in \R^n$ in terms of its $1$-norm:
$$
     |x| := \sum_{j=1}^n |x_j|.
$$

Our main result is:

\medskip
\begin{theorem}\label{t.main} The Ford-Fulkerson method for the max-flow problem (\ref{eq.flow-problem}) in a regular space $\cV$ requires at most $|E|^2$ augmentations if always augmenting paths $P\in \cP$ of shortest possible length  $|P|$ are selected.
\end{theorem}

For the proof of Theorem~\ref{t.main}, we need to establish some more structural properties of the path system $\cP$  in a regular space.

\subsection{Path algebra}\label{sec:path-algebra}
Let $P,Q\in \cP$ be arbitrary $r$-paths.  By Lemma~\ref{lemma2}, $P+Q$ is a sum of primitive vectors $p^i$, each conforming to $P+Q$:
$$
P+Q = p^1 +p^2 +\ldots +p^k.
$$
In view of $(P+Q)_r = P_r+Q_r = 2$, exactly two of these summands $p^i$, say $P\wedge Q$ and $P\vee Q$, have $(P\wedge Q)_r = +1 = (P\vee Q)_r$ and thus are $r$-paths. By conformity,  we moreover find for each $j\in E$:
$$
\begin{array}{cccl}
(a) &(P\wedge Q)_j \neq 0  &\Longrightarrow& (P\wedge Q)_j = P_j \;\mbox{or}\;
(P\wedge Q)_j = Q_j.\\
(b) &(P\vee Q)_j \neq 0 &\Longrightarrow& (P\wedge Q)_j = P_j \;\mbox{or}\;
(P\wedge Q)_j = Q_j.\\
(c) &P_j Q_j =-1 &\Longrightarrow& (P\wedge Q)_j = (P\vee Q)_j = 0.\\
(d) &(P\wedge Q)_j\cdot(P\vee Q)_j \neq 0 &\Longrightarrow& (P\wedge Q)_j=(P\vee Q)_j =P_j = Q_j.
\end{array}
$$
We say that $P$ and $Q$ are \emph{nonconformal} if $P_j Q_j=-1$ holds for some $j\in E$. Otherwise, $P$ and $Q$ are \emph{conformal}.

\medskip
\begin{lemma}\label{l.length-inequality} Let $P,Q\in \cP$ be arbitrary paths. Then
\begin{itemize}
\item[(i)] $ |P+Q| = |P| + |Q|\quad \Longleftrightarrow\quad$ $P$ and $Q$ are conformal.
\item[(ii)] $ |P\wedge Q| +|P\vee Q| \leq  |P| +|Q|$.
\item[(iii)] $|P\wedge Q| +|P\vee Q| < |P| +|Q|$ if $P$ and $Q$ are nonconformal.
\end{itemize}
\end{lemma}

{\small\medskip\Pf (i) is an obvious property of the $1$-norm.
To see (ii), assume
$$
P+Q = p^1 +\ldots +p^k,
$$
where the $p^i$ are primitive vectors each conforming to $P+Q$. By the triangle inequality of the $1$-norm, we therefore find
$$
|P\wedge Q| +|P\vee Q|  \leq  |p^1|+\ldots + |p^k| =|P+Q| \leq |P|+|Q|.
$$
By (i), the latter inequality is strict if $P$ and $Q$ are nonconformal.

\qed}

\subsection{Shortest augmenting paths}\label{sec.shortest-augmenting-paths} Let $f\in \cV$ be a feasible solution for the max-flow problem (\ref{eq.flow-problem}). Let $\cP^f\subseteq \cP$ the collecting of $f$-augmenting $r$-paths $P$ of shortest possible length $|P|$.

\medskip
\begin{lemma}\label{l.shortest-conformal1} For any $P,Q\in \cP^f$, one has $P\wedge Q\in \cP^f$ and $P\vee Q\in \cP^f$.
\end{lemma}

{\small\medskip Suppose $(P\vee Q)\notin \cP^f$, for example. Then there is some $j\in E$ such that $f_j = 0$ and $(P\wedge Q)_j = -1$ or $f_j = c_j$ and $(P\wedge Q)_j = +1$. Without loss of generality, let us assume the former (the latter can be dealt with analogously). So property (a) above implies $P_j =-1$ or $Q_j = -1$, which contradicts the hypothesis that $P$ and $Q$ are augmenting.

\medskip
If both $P\wedge Q$ and $P\vee Q$ are augmenting, we have $|P\vee Q| \geq |P|$ and $|P\vee Q|\geq |Q|$ and thus conclude from
$$
     |P| + |Q| \leq  |P\wedge Q|+|P\vee Q| \leq |P|+|Q|
$$
that equality must hold, which means that $P\wedge Q$ and $P\vee Q$ are shortest augmenting paths as well.

\qed}

\medskip
\begin{lemma}\label{l.shortest-conformal2} Any two paths $P,Q\in \cP^f$ are conformal.
\end{lemma}

{\small\medskip\Pf By Lemma~\ref{l.shortest-conformal1}, we know $P\wedge Q, P\vee Q\in \cP^f$.  Hence, if $P$ and $Q$ were nonconformal, we would arrive at the contradicition
$$
  |P| +|Q| = |P\wedge Q|+|P\vee Q| < |P|+|Q|.
$$
\qed}

\medskip
\begin{lemma}\label{l.shortest-augmentation} Let $f$ be a feasible flow and $g$ the feasible flow obtained by augmenting $f$ along a shortest path $P\in \cP^f$. Let $Q\in \cP^g$ be a shortest augmenting path for $g$. Then
$$
|Q| > |P| \quad \Longleftrightarrow \quad Q\notin \cP^f.
$$
\end{lemma}

{\small\medskip\Pf If $Q\in \cP^f$, we have $|Q| =|P|$ by the definition of $\cP^f$. Consider  now the case $Q\in \cP^g\setminus\cP^f$. If $Q$ is augmenting for $f$, then $|Q|>|P|$ follows from the definition of $\cP^f$. It remains to analyze the case where $Q$ is not augmenting for $f$.

\medskip
We claim $P\wedge Q \in \cP^f$ and $P\vee Q\in \cP^f$. Suppose $P\wedge Q\notin \cP^f$, for example. So some $j\in E$ exists with either $(P\wedge Q)_j = -1$ and $f_j = 0$ or $(P\wedge Q)_j =+1$ and $f_j =c_j$. Let us assume the latter (the former can be dealt with in the same way). Then $P_j \leq 0$ must hold (otherwise, $P$ would not be augmenting for $f$). Because $P\wedge Q$ conforms to $P+Q$, we therefore conclude
$$
 (P+Q)_j= +1 \quad\mbox{and hence}\quad    Q_j =+1 \;\mbox{and}\; P_j = 0.
$$
But this implies $g_j = f_j$ and thus $Q\notin \cP^g$, contradicting our hypothesis. So the claim is established.

\medskip
On the other hand, since $Q$ augments $g$ but not $f$, $P$ and $Q$ must be conformal. Lemma~\ref{l.length-inequality} therefore yields
$$
|P|+|P| = |P\wedge Q| + |P\vee Q| < |P| + |Q| \quad\mbox{and thus}\quad |P| < |Q|.
$$
\qed}

\paragraph{Proof of Theorem~\ref{t.main}.}  We are now in the position to establish Theorem~\ref{t.main}. If we augment along a shortest path $P$, the next shortest augmenting path $Q$ will

\begin{itemize}
\item[(i)] either be of the same length $|Q|=|P|$  and(!) conformal with $P$
\item[(ii)]  or be of length $|Q|\geq |P|+1$.
\end{itemize}

\medskip
On the other hand, augmentation along paths that are conformal with the previous one cannot happen more than $|E|$ times in a row because after each conformal augmentation the flow becomes tight on some $j\in E$ at its lower capacity bound $0$ or its upper bound $c_j$ and stays tight until a nonconformal augmentation occurs.

\medskip
Since all paths $P\in \cP$ have length $|P|\leq |E|$, no more than $|E|-1$ nonconformal augmentations will take place, which yields the bound $|E|^2$ on the total number of augmentations.

\section*{Final remarks} Let $G=(V,E)$ be a directed graph with distinguished arc $r=(t,s)$ and vertex-arc incidence matrix $A$. The support $\|p\|$ of a primitive vector $p\in \ker A$  is a circuit in the circuit matroid $\cM(G)$ associated with $G$. An $r$-path relative to $\cV=\ker A$ is thus a circuit with oriented arcs that contains $r$ and is traversed so that $r$ is a ''forward arc''. A shortest $r$-path thus corresponds to a shortest directed path from the vertex $s$ to the vertex $t$ in the auxiliary graph $G(f)$ of forward and backward arcs relative to a given feasible flow $f$.

\medskip
A circuit of $G$ contains at most $|V|$ arcs. So the Edmonds-Karp bound $|V|\cdot|E|$ on the number of shortest path augmentations follows for the classical Ford-Fulkerson algorithm.

\medskip
In the case of co-flows on $G$, our regular space $\cV$ is the row space of $A$. The supports $\|p\|$ of the primitive vectors $p\in \cV$ are then the minimal cutsets of $G$ (\emph{i.e.}, the co-circuits of the circuit matroid $\cM(G)$ of $G$). An $r$-path is a signed cutset containing $r$. A shortest augmenting $r$-path relative to a given feasible co-flow $f$ can thus be computed as a minimal $r$-cut in the associated auxiliary graph $G(f)$ of forward and backward arcs relative to $f$.

\end{document}